\documentclass[a4paper,11pt]{elsarticle} 
\usepackage{lineno,hyperref}
\modulolinenumbers[5]
\usepackage{geometry}
\usepackage{hyperref}
\usepackage{amsfonts}
\usepackage{amsmath}
\usepackage{amsthm}
\usepackage{algorithm}  
\usepackage{algpseudocode}  
\usepackage{amsmath}  
\usepackage{graphicx}
\usepackage{subfigure}
\usepackage{booktabs}  
\usepackage{multicol}  
\usepackage{multirow}
\usepackage{threeparttable}  
\usepackage{makecell}
\usepackage{caption2}
\usepackage{lipsum}
\usepackage{bm}
\usepackage{multirow}
\usepackage{amssymb}

\date{}

\newtheorem{definition}{Definition}

\usepackage{color}

\def\tcb{\textcolor[rgb]{0,0,1}}

\usepackage{url}

\makeatletter
\def\ps@pprintTitle{%
  \let\@oddhead\@empty
  \let\@evenhead\@empty
  \let\@oddfoot\@empty
  \let\@evenfoot\@oddfoot
}
\makeatother

\begin{document}
	
\begin{frontmatter}
\title{A Biased Deep Tensor Factorization Network For Tensor Completion\tnoteref{mytitlenote}}
\tnotetext[mytitlenote]{Research supported by the National Natural Science Foundation of China under Grant (11801418).}

\author[1]{Qianxi Wu}
\ead{wuqx214@gmail.com}
\address[1]{College of Computer and Artificial Intelligence, Wenzhou University, Zhejiang 325035, China.}

\author[2]{An-Bao Xu\corref{mycorrespondingauthor}}
\cortext[mycorrespondingauthor]{Corresponding author}
\ead{xuanbao@wzu.edu.cn}

\address[2]{College of Mathematics and Physics, Wenzhou University, Zhejiang 325035, China.}

\begin{abstract}
	Tensor decomposition is a popular technique for tensor completion, However most of the existing methods are based on linear or shallow model, when the data tensor becomes large and the observation data is very small, it is prone to over fitting and the performance decreases significantly. To address this  problem, the completion method for a tensor based on a Biased Deep Tensor Factorization Network (BDTFN) is proposed. This method can not only overcome the shortcomings of traditional tensor factorization, but also deal with complex non-linear data. Firstly, the horizontal and lateral tensors corresponding to the observed values of the input tensors are used as inputs and projected to obtain their horizontal (lateral) potential feature tensors. Secondly, the horizontal (lateral) potential feature tensors are respectively constructed into a multilayer perceptron network. Finally, the horizontal and lateral output tensors are fused by constructing a bilinear pooling layer. Tensor forward-propagation is composed of those three step, and its parameters are updated by tensor back-propagation using the multivariable chain rule. In this paper, we consider the large-scale 5-minute traffic speed data set and use it to address the missing data imputation problem for large-scale spatiotemporal traffic data. In addition, we compare the numerical performance of the proposed algorithm with those for state-of-the-art approaches on video recovery and color image recovery.  Numerical experimental results illustrate that our approach is not only much more accurate than those state-of-the-art methods, but it also has high speed.\\
\end{abstract}

 
\end{frontmatter}

\section{\large\bf Introduction}

In the field of science and engineering, data defect is common. The goal of tensor completion is to predict the missing elements by learning some known data. Tensor analysis has reached a celebrity status in the areas of machine learning, computer vision, and artificial intelligence \cite{WuTanLi2017}. Tensor completion is an important problem for tensor analysis. Many problems can be solved by tensor completion method, especially in time series prediction \cite{ChenXY2020}. Other popular applications include computer vision \cite{LMWY2013} \cite{X2020} \cite{LTYL2016}, dimensionality reduction \cite{ZJ2016}, signal processing \cite{ZEAHK2014} and data mining \cite{KS2008} \cite{SPLCLQ2009}. With the significant progress of low-cost sensing technology, more and more real-world traffic measurement data acquisition has high space sensor coverage and good time resolution \cite{ChenXY2020}. These large-scale and high-dimensional traffic data provide us with unprecedented opportunities to perceive traffic dynamics and develop efficient and reliable applications for its. One of the challenges is the loss and damage of data, which is not conducive to obtaining real signals. The other is the computational challenge of using existing methods to analyze large-scale and high-dimensional traffic data. The problem of data loss and damage may result from complex sensor failure, communication failure, or even low sensor coverage. In order to improve data quality and support downstream applications, missing data filling becomes an important task.  Despite recent efforts and advances in developing machine learning-based imputation models, most studies still focus on small-scale data sets and it remains a critical challenge to perform accurate and efficient imputation on large-scale traffic data sets \cite{ChenXY2020}.

According to the literature \cite{ChenXY2020}, there are numerous methods for the missing data imputation problem in spatiotemporal traffic data. By using tensor structure to model multivariable time series data, multivariable / multidimensional settings are studied. Signature framework includes: 1) classical time series analysis methods, such as vector autoregressive model \cite{Bashir2018} , 2) pure low rank matrix factorization / completion, such as low rank matrix factorization \cite{Asif2013,Asif2016}, principal component analysis \cite{Qu2008,Qu2009} and its variants, and 3) Pure low rank tensor factorization / completion, such as Bayesian tensor factorization \cite{Chen2019} and low rank tensor factorization \cite{Ran2016,ChenYang2020} , 4) by integrating low rank matrix factorization of time series models, such as time regularized matrix factorization \cite{Yu2016} and Bayesian time matrix factorization \cite{Sun2019}, and 5) low rank tensor completion by integrating time series models (such as LATC) \cite{ChenSun2020}. For this large class of methods, a common goal is to capture temporal dynamics and spatial consistency. For example, by integrating the time series model and the pure low rank tensor completion model, the recently developed LATC framework shows superior performance compared with the existing low rank tensor completion model \cite{ChenSun2020} by imposing local consistency and global consistency on the traffic series data \cite{LiSu2015}. In addition, some nonlinear methods (such as deep learning in \cite{ZHZZJ2021,CL2019,Che2018}) have been applied to solve the problem of missing traffic data. 

In the tensor completion algorithm, the most commonly used method is tensor decomposition and its improved method, such as reference \cite{Yongming2016,ZA2017,X2020,LMWY2013}. However, these models assume that the defective data tensor is of low rank and can be generated linearly by latent eigenvectors, that is, the prediction of defective data is directly obtained by the linear interaction of latent eigenvectors. However, when the missing data become complex and diverse, such as the known data set is not a certain random distribution or the number of known observation data is too small \cite{Candes2009}, the effect is often not ideal. In recent years, with the improvement of computer computing ability, the increasing of data quantity and the breakthrough of deep learning in image classification field, many scholars apply the existing deep learning methods to the tensor completion, such as FCTN-TC \cite{ZHZZJ2021}, LATC-Tubal \cite{ChenXY2020} and ATF \cite{CL2019}, and obtained better results than shallow model. At present, the tensor completion algorithms based on deep learning generally includes multiple network layers, such as input layer, embedding layer, non-linear transformation layer (such as Sigmoid and ReLU), full connection layer and output layer, etc. Different from shallow tensor completion, deep learning based tensor completion generally gets the predicted value of the model by multi-layer nonlinear mapping of potential eigenvectors. In this paper, the deep learning technique is applied to tensor decomposition. Generally, the horizontal and lateral corresponding elements of data tensor are embedded into the latent matrix space respectively. Since the data are usually located on a low dimensional manifold, the dimension of the embedding tensor is generally much smaller than the order of the data tensor. In shallow tensor completion, latent tensor is used to construct linear model directly. Because of the direct linear interaction of latent tensor, these methods need to attach many conditions to the feature space of the known data in order to get better recovery results. As mentioned in reference \cite{Candes2010}, there are requirements for the number and probability distribution of the known observation data. Obviously, it is difficult to meet these conditions in practical application. In this regard, reference \cite{SiS2016} proposed nonlinear transformation of potential variables to reduce these limitations and include side information into the model. Deep learning theory is one of the key technologies in the field of artificial intelligence  \cite{GBC2016}. Because it can learn the complex nonlinear structural features of data hiding. The tensor completion algorithm based on deep learning usually transforms the latent tensor nonlinearly and trains it by constructing a multilayer neural network. The input of the network is usually the latent feature tensor, and the output is the model prediction value of the element. For example, in the depth tensor completion algorithm proposed in reference \cite{Fan2017}, both the network parameters and the input latent vectors are learned as unknown parameters. The loss function is constructed by the actual values of known elements and the values obtained from the model. The hidden variables and the network parameters are optimized at the same time, and the numerical solutions of the potential variables and the network parameters are obtained by the random gradient descent method.

In this paper, a Biased Deep Tensor Factorization Network (BDTFN) is proposed based on deep learning and tensor decomposition, inspired by the idea of a Biased Deep Matrix Factorization (BDMF) model  \cite{MLPLD2019}. The method uses multi-layer perceptron to extend the latent eigenvectors to multi-layer neural network.
As an improvement of reference \cite{XueHJ2017,Fan2017,Dziugaite2015}, this paper mainly has the following four innovations:

\begin{itemize}

\item A multi-layer perceptron network is constructed from the horizontal (lateral) potential feature tensors which is obtained by projecting The horizontal (lateral) tensors of the input tensors.

\item A bilinear pooling layer is designed to use Hadamard product to fuse the horizontal and lateral network outputs, and then obtain the predicted model value.

\item A tensor back-propagation is given to update our parameters by computing the gradients using the multivariable chain rule.

\item Through the combining of tensor forward-propagation and tensor back-propagation, a new tensor completion algorithm is proposed, named deep tensor factorization network. This method extends traditional tensor decomposition to multi-layer neural network, and then have ability to deal with nonlinear data with complex structure.

\end{itemize}

\section{\large\bf Background and preliminaries}



In this paper, we denote vectors by lowercase boldface letters, $\mathbf{a}$, and matrices by uppercase boldface letters, $\mathbf{A}$, respectively. Tensors are represented in boldface Euler script letters. For instance, $\bm{\mathcal{A}} \in \mathbb{R} ^{n_{1}\times n_{2}\times n_{3}   } $ is used to denote a third-order tensor, and its $\left ( i,j,k \right ) $th entry is represented as $\bm{\mathcal{A}} _{i,j,k} $ or $a _{i,j,k} $. The $i$th horizontal, lateral and frontal slice (see definitions in \cite{KB2009}) of third-order tensor $\bm{\mathcal{A}}$ is denoted by the Matlab notation $\bm{\mathcal{A}}\left ( i,:,: \right ) $, $\bm{\mathcal{A}}\left ( :,i,: \right ) $ and $\bm{\mathcal{A}}\left ( :,:,i \right ) $ respectively. In general, $\bm{\vec{\mathcal{A}}}_j  $ denotes the lateral slice $\bm{\mathcal{A}}\left ( :,j,: \right )$ compactly, and $\bm{\mathcal{A}}\left ( i,j,: \right ) $ denotes the tube of $i,j$ in the third-order tensor dimension.


By avoiding vectorization of the data, our BDTFN frame work has the potential to extract multidimensional correlations. We introduce our design for third-order tensors (i.e., three-dimensional arrays), but the theory easily extends to higher dimensions \cite{CDMartin2013,NHao2014}. We will demonstrate that using the $*_{M}$ -product framework creates a more efficient, powerful parameterization than its matrix equivalent and its matrix-mimeticity enables an elegant high-dimensional extension. 

\begin{definition}\label{Mode-3 product} 
	\textbf{(Mode-3 product) \cite{N2019}} Suppose $\bm{\mathcal{A}}\in \mathbb{R}^{n_{1}\times n_{2}\times n_{3}}$ and $\mathcal{M}\in \mathbb{R}^{p\times n_{3}}$. Then the mode-3 product is equivalent to the following:
	\begin{center}
	$	\bm{\mathcal{A}}\times_{3}\bm{\mathcal{M}} \equiv (\bm{\mathcal{M}}\otimes\bm{\mathcal{I}}_{n_{1}})\cdot\bm{\mathcal{A}}_{(2)}$
	\end{center}
	where $\otimes$ denotes the Kronecker product and $\bm{I}_{n_{1}}$ is the $n_{1}\times n_{1}$ identity matrix and $\bm{\mathcal{A}}_{(2)}$ is the Mode-2 unfolding.
	
\end{definition}
\begin{definition}\label{Facewise product}
	\textbf{(Facewise product) \cite{N2019}} Given $\bm{\mathcal{A}}\in \mathbb{R}^{l\times p\times n}$ and $\bm{\mathcal{B}}\in \mathbb{R}^{p\times m \times n}$. Then the facewise product multiplies the frontal slices of $\bm{\mathcal{A}}$ and $\bm{\mathcal{B}}$, denoted by the \tcb{$ "\vartriangle"$} operation, as follows:
	\begin{center}
		$\bm{\mathcal{C}}=\bm{\mathcal{A} \vartriangle
			\mathcal{B}}\equiv\bm{\mathcal{C}}^{(i)}
		=\bm{\mathcal{A}}^{(i)}\cdot\bm{\mathcal{B}}^{(i)}\quad\quad for\quad i=1,...,n.$
	\end{center}	
\end{definition}
\begin{definition}\label{*_{M}-product}
	\textbf{($*_{M}$-product) \cite{N2019}} Given $\bm{\mathcal{A}}\in \mathbb{R}^{l\times p\times n}$ and $\bm{\mathcal{B}}\in \mathbb{R}^{p\times m \times n}$ and an invertible $n\times n$ matrix $\bm{\mathcal{M}}$, then the $*_{M}$-product is defined as follows :
	\begin{center}
		$	\bm{\mathcal{C}}=\bm{\hat{\mathcal{A}}} *_{M}
		\bm{\hat{\mathcal{B}}}=(\bm{\mathcal{A}\vartriangle \mathcal{B}})\times_{3}\bm{\mathcal{M}}^{-1}, \quad \bm{\mathcal{C}}\in \mathbb{R}^{l\times m\times n},$
	\end{center}
	where $\bm{\hat{\mathcal{A}}}=\bm{\mathcal{A}}\times_{3}\bm{\mathcal{M}}$ and $\bm{\hat{\mathcal{B}}}=\bm{\mathcal{B}}\times_{3}\bm{\mathcal{M}}$ are the tensors in the transform domain.
\end{definition}
	Let $\bm{\mathcal{M}}=\bm{\mathcal{F}}_{n}$ where $\bm{\mathcal{F}}_{n}$ is the unmornalized $n \times n$ discrete Fourier transform(DFT) matrix. 
	Then the $*_{M}$-product is equivalent to the t-product \cite{MEKCDM2011}.

After we evaluate the performance of our network, we update our parameters by computing the gradients using the multivariable chain rule (5.5) of \cite{N2019}. We will derive the tensor backward propagation formulas to further demonstrate the elegance of working in a matrix-mimetic framework. We first derive some basic $*_{M}$ -product gradients using \cite{PP2012} as reference for matrix calculus. We start with the $*_{M}$-product gradient.

Based on Definition 5.2.5 and Definition 5.2.6 in \cite{N2019}, the back-propagation formula for the $*_{M}$-product is proposed as follows,.
\begin{definition}\label{*_{M}-product gradient} 
	\textbf{($*_{M}$-product gradient) \cite{N2019}} Suppose $f: \mathbb{R}^{l\times m\times n} \rightarrow \mathbb{R}$ is a scalar-valued, differentiable function and let $\bm{M}$ is an invertible $n\times n$ matrix. Let $\bm{\mathcal{A}} \in \mathbb{R}^{l\times p\times n}$ and $\bm{\mathcal{B}} \in \mathbb{R}^{p\times m\times n}$. Then, 		
\begin{equation}
	\begin{aligned}		
	\dfrac{\partial}{\partial \bm{\mathcal{B}}} [f(\bm{\mathcal{A}} *_{M} \bm{B})]  &= \dfrac{\partial}{\partial \vec{\bm{\mathcal{B}}}}[f(((\bm{\mathcal{A}} \times_{3} \bm{M}) \vartriangle (\bm{\mathcal{B}} \times_{3} \bm{M})) \times_{3} \bm{M}^{-1})] \\
  &=[(\bm{\mathcal{A}}\times_{3} \bm{M})^{\intercal} \vartriangle (f'(\bm{\mathcal{A}} *_{M} \bm{\mathcal{B}}) \times_{3} \bm{M}^{\intercal})]\times_{3} \bm{M}^{\intercal}.
  \end{aligned}
\end{equation}	
\end{definition}

To simplify the formula above, consider the case when $\bm{M}$ is a nonzero multiple of an orthogonal matrix so that $\bm{M}^{-1}=\bm{M}^{\intercal}$. Then,	
\begin{equation}
\dfrac{\partial}{\partial \bm{\mathcal{B}}} [f(\bm{\mathcal{A}} *_{M} \bm{\mathcal{B}})] =\bm{\mathcal{A}}^{\intercal} *_{M} f'(\bm{\mathcal{A}} *_{M} \bm{\mathcal{B}}).		
\end{equation}

\section{\large \bf Biased Deep Tensor Factorization Network}

\begin{figure}[H]
	\centering
	\noindent\makebox[\textwidth][c] {
		\includegraphics[scale=0.25]{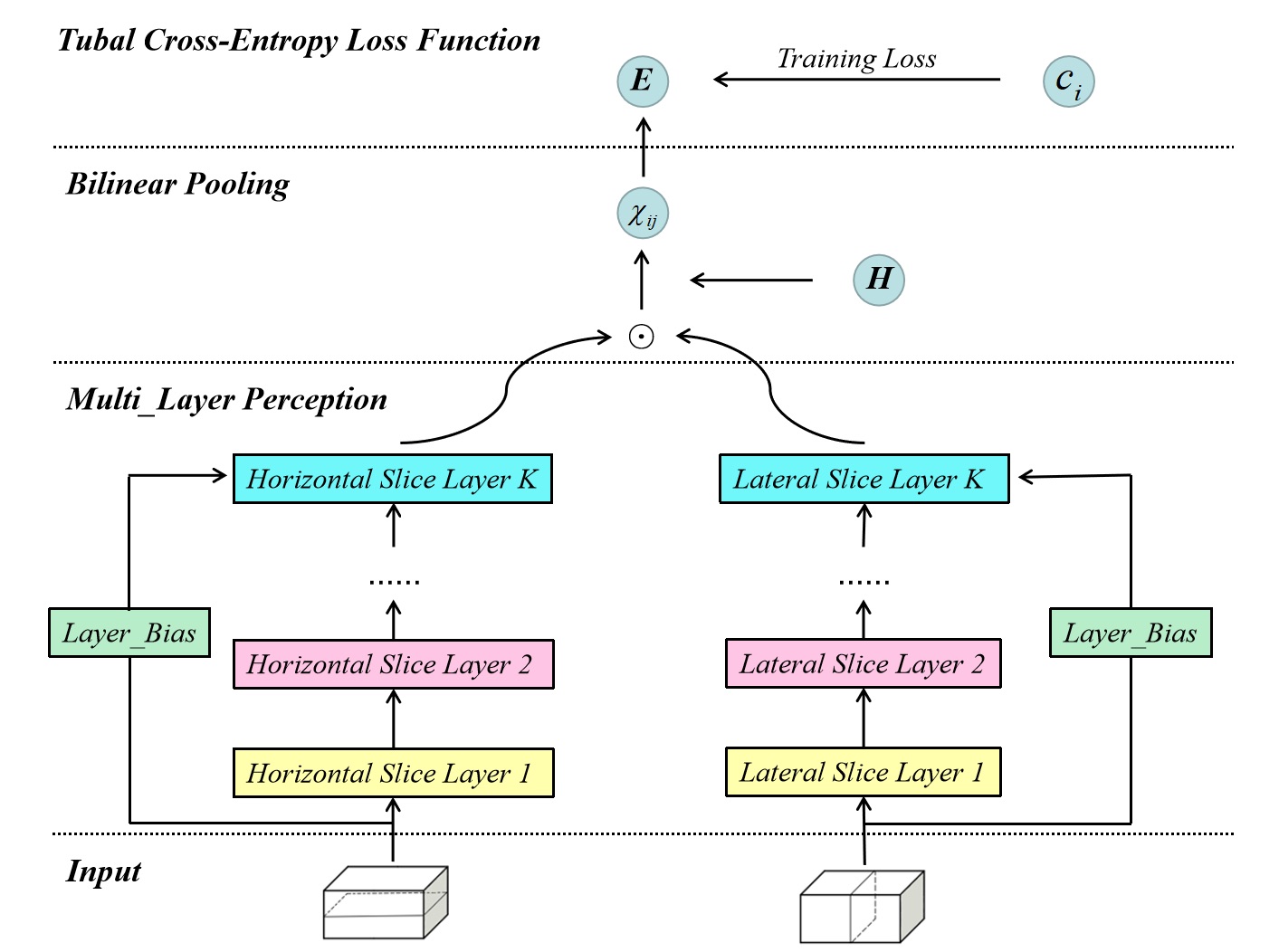} }
	\caption{The Architecture of Biased Deep Tensor Factorization Network Model.}
	\label{DTFN Model}
\end{figure}

The BDTFN model improves the nonlinear processing ability by combining different layers. The algorithm proposed in this paper is shown in Figure \ref{DTFN Model}, which mainly includes input layer,  multi-layer perceptron layer, bilinear pooling layer and output layer. In the training stage, each observation point is transformed by its corresponding horizontal and lateral tensors to obtain the potential feature tensors, which are then used as the input of the multi-layer perceptron layer, and fused by the horizontal and lateral potential feature tensors of the bilinear pooling layer. Finally, the square loss function is constructed by using the model value and the known observation value. Tensor forward-propagation is composed of those three step, and its parameters are updated by tensor back-propagation using the multivariable chain rule. 

\subsection{Tensor Forward-Propagation}
\subsubsection{Input layer}
Firstly, the complete data in the data set is processed to get the tensor $\bm{\mathcal{R}} \in \mathbb{R}^{n_{1}\times n_{2}\times n_{3}}$. Secondly, the processed tensor is transformed horizontally (\bm{$\mathcal{U}$}) and laterally (\bm{$\mathcal{V}$}). The two transformed tensors are used as inputs.

Latent Horizontal Slice
\begin{equation}\label{Horizontal}
	\bm{\vec{\mathcal{U}}_{i}^{1}} = \bm{\mathcal{R}}(i,:,:)^{T},
\end{equation}

Latent Lateral Slice
\begin{equation}\label{Lateral}
	\bm{\vec{\mathcal{V}}_{i}^{1}} = \bm{\mathcal{R}}(:,i,:),
\end{equation}
where $\bm{\vec{\mathcal{U}}_{i}^{0}}  \in \mathbb{R}^{n_{2}\times 1 \times n_{3}}$ and $\bm{\vec{\mathcal{V}}}_{i}^{0}  \in \mathbb{R}^{n_{1}\times 1 \times n_{3}}$ are the horizontal and lateral potential characteristic tensors obtained by nonlinear transformation.

\subsubsection{Multi-layer Perception}
Multi-Layer perceptron is mainly used to make multi-layer nonlinear mapping of potential eigenvalues, so as to better learn the nonlinear characteristics of data. In this layer, the potential eigenvalues are fully connected by multi-layer perceptron.

Latent Horizontal Slice

\begin{equation}
	\bm{\vec{\mathcal{U}}_{i}^{j+1}} = \sigma(\bm{\mathcal{W}_{u}^{j}} {\ast}_{M}  \bm{\vec{\mathcal{U}}_{i}^{j}} + \bm{\vec{\mathcal{B}}_{u}^{j}}),j=1,...,N-1
\end{equation}

\begin{equation}\label{key4}
	\bm{\vec{\mathcal{U}}_{i}^{N}} = \sigma(\bm{\mathcal{W}_{u}^{j}} {\ast }_{M}  \bm{\vec{\mathcal{U}}_{i}^{N-1}})
\end{equation}
where $\sigma$ is an element-wise, momlinear activation function and $N$ is the number of layers in the network. Note that $\bm{\mathcal{W}_{u}^{j}} \in \mathbb{R}^{ m \times n_2 \times n_3}$ and is a weight tensor and $\bm{\vec{\mathcal{B}}_{u}^{j}} \in \mathbb{R}^{m \times1\times n_3}$ is a bias lateral slice.

Latent Lateral Slice

\begin{equation}
	\bm{\vec{\mathcal{V}}_{i}^{j+1}} = \sigma(\bm{\mathcal{W}_{v}^{j}} {\ast}_{M}  \bm{\vec{\mathcal{V}}_{i}^{j}} + \bm{\vec{\mathcal{B}}_{v}^{j}}),j=1,...,N-1
\end{equation}

\begin{equation}\label{key6}
	\bm{\vec{\mathcal{V}}_{i}^{N}} = \sigma(\bm{\mathcal{W}_{v}^{j}} {\ast}_{M}  \bm{\vec{\mathcal{V}}_{i}^{N-1}})
\end{equation}
where $\bm{\mathcal{W}_{v}^{j}} \in \mathbb{R}^{ n \times n_2 \times n_3}$ is a weight tensor and $\bm{\vec{\mathcal{B}}_{v}^{j}} \in \mathbb{R}^{n \times1\times n_3}$ is a bias lateral slice.

Among them, \bm{$\mathcal{W}_{u}$} and \bm{$\mathcal{W}_{v}$} are the horizontal and lateral characteristic projection tensors of the tensor, which are unknown training parameters and need to be trained by back propagation; $\sigma(\cdot)$ is a nonlinear activation function, which mainly adds nonlinear factors to the linear transformation to make up for the deficiency of the linear model, which is the difference between the deep decomposition network and the traditional tensor decomposition, and also an important factor for the success of deep learning. The common activation functions are Sigmoid, ReLU, tanh and so on. However, in tensor decomposition, the effect of Sigmoid is usually better than that of other activation functions. The experiment in reference \cite{Sedhain2015} also has a similar conclusion, which will be applied in the experimental part of this paper. The activation function does not change the tensor dimension.
Because the multi-layer perceptron is essentially a fully connected input, and this paper defines it as a non-linear mapping with equal dimensions, the number of layers should not be too high, otherwise it is easy to over fit.

\subsubsection{Bilinear Pooling}

After multi-layer perceptron mapping of the potential eigenvalues of the horizontal and lateral tensors, the next step is to fuse the horizontal and lateral network outputs to get the predicted model value, using the $*_{M}$-product. In order to further improve the generalization and learning ability, this paper constructs bilinear pooling layer.

\begin{equation}
	{\bm{x}}_{ij} = \bm{\mathcal{X}}(i,j,:) =  \bm{\vec{\mathcal{H}}^{T}} {\ast }_{M} \sigma(\bm{\vec{\mathcal{U}}_{i}^{N}} \odot \bm{\vec{\mathcal{V}}_{j}^{N}}),
\end{equation}
where the notation $\odot$ represents Hadamard product, and $\mathcal{H}$ is the unknown parameter to be trained.

Or use the following bilinear pooling layer,
\begin{equation}\label{key7}
\bm{\mathcal{X}} = (\bm{\mathcal{U}}^N)^T  *_M \bm{\mathcal{V}}^N.
\end{equation}

\subsubsection{Tubal Cross-Entropy Loss Function}

Using the tubal softmax function ${ [ f(\bm{\vec{\mathcal{X}}}) ] }_{i}$ of a lateral slice$\bm{\vec{\mathcal{X}}}$ in \cite{N2019} , we have
\begin{equation}
	\bm{\vec{\mathcal{Y}}} \equiv { [ f(\bm{\vec{\mathcal{X}}}) ] }_{i} = (\sum\limits_{j=1}^{n_{1}} \exp( \bm{x}_{j}))^{-1} {\ast}_{M} \exp( \bm{x}_{i}),
\end{equation}
where $\bm{\vec{\mathcal{Y}}} \in \in \mathbb{R}^{ n_1 \times 1 \times n_3}$

Then we obtain
\begin{equation}
	H(\bm{\vec{\mathcal{Y}}}_{i}, c_{i}) = -\log\bm{y_{i}}, c_{i}
\end{equation}
where $H( \cdot, \cdot)$ is a tubal cross-entropy function and $c$ is the true classification label (i.e. the true distribution is Bernoulli) in the reference \cite{N2019}.

So we get our tubal cross-entropy loss function as follows:
\begin{equation}\label{key10}
	E(\bm{\mathcal{N}}) = \dfrac{1}{2n_{2}}\sum\limits_{i=1}^{n_{2}} \lVert	H(\bm{\vec{\mathcal{Y}}}_{i}, c_{i})\rVert_{F}^{2} + \mathcal{R}(\bm{\mathcal{N}})
\end{equation}
where \bm{$x$},\bm{$y$} $\in$ $\mathbb{R}^{1\times1\times{n_{3}}}$ are the tubes. Let $E$ is our objective function.

\subsection{Tensor Back-Propagation (Multivariable Chain Rule)}

In this subsection, we determine the back-propagation formula by computing the derivatives for the tensor forward-propagation. The most notable aspect of these derivations is that tubes will continue to mimic scalars in our $*_{M}$-framework; i.e., our derivatives will look nearly identical to derivatives of one-dimensional functions.

\subsubsection{Derivation for Tubal Cross-Entropy Loss Function}

We will give the back-propagation formula for our tubal cross-entropy loss function. For formulaic simplicity, we assume the transformation matrix $\bm{M}$ is a non-zero multiple of an orthogonal matrix. 

We will derive the derivative of the tubal cross-entropy loss function using the following derivations in \cite{N2019}.

\begin{definition}\label{Tubal cross-entropy}
	\textbf{(Tubal cross-entropy) \cite{N2019}} Let $\bm{y}_{i}=(\sum_{j=1}^{p}\exp(\bm{x}_{j}))^{-1} *_{M} \exp(\bm{x}_{i})$ be the output from the tubal softmax function. Then,
	\begin{equation}
		\begin{aligned}		
			\dfrac{\partial}{\partial\bm{x}_{k}}(\dfrac{1}{2}\lVert-\log(\bm{y}_{i})\rVert_{F}^{2}) 
			&=\dfrac{\partial}{\partial\bm{x}_{k}}(-\bm{y}_{i}^{-1})*_{M}(-\log(\bm{y}_{i}))
		\\	&=\begin{cases}
				-\bm{y}_{k} *_{M} \log(\bm{y}_{i}) \qquad\quad i \neq k \\
				-(\bm{y}_{k}-\bm{e}) *_{M} \log(\bm{y}_{k})\quad i = k,
			\end{cases}
		\end{aligned}
	\end{equation}
\end{definition}
Thus, the derivative of our tubal cross-entropy loss function is the following:


\begin{equation}\label{key_delta_x}
	\begin{aligned}
	\delta \bm{\mathcal{X}} = \dfrac{\partial E}{\partial \bm{\mathcal{X}}} \equiv \dfrac{\partial}{\partial x_{k}}(\dfrac{1}{2} \lVert -\log(\bm{y_{i}}) \rVert_{F}^{2})
	= \dfrac{\partial}{\partial x_{k}}(\bm{y_{i}}) {\ast}_{M} (\bm{-y_{i}^{-1}}) {\ast}_{M} (-\log(\bm{y_{i}}))  \\
	 = \begin{cases}
		\bm{-y_{k}} {\ast}_{M} \log(\bm{y_{i}}), \quad\quad\quad i \neq k \\
		\bm{(-{y_{k}}-e)} {\ast}_{M} \log(\bm{y_{i}}), \quad i=k
	\end{cases}
	\end{aligned}
\end{equation}

\subsubsection{Derivation for Bilinear Pooling}

Derivation formulas for bilinear pooling are proposed by using Definition 4 as follows.

\begin{equation}\label{key_delta_U^N}
	\delta \bm{\mathcal{U}}^{N} = \dfrac{\partial E}{\partial \bm{\mathcal{U}}^{N}} = \dfrac{\partial E}{\partial \bm{\mathcal{X}}^{T}} {\ast}_{M} \dfrac{\partial \bm{\mathcal{X}}^{T}}{\partial \bm{\mathcal{U}}^{N}} = \bm{\mathcal{V}}^{N} {\ast}_{M} (\delta \bm{\mathcal{X}})^{T}
\end{equation}

\begin{equation}\label{key_delta_V^N}
	\delta \bm{\mathcal{V}}^{N} = \dfrac{\partial E}{\partial \bm{\mathcal{V}}^{N}} = \dfrac{\partial E}{\partial \bm{\mathcal{X}}^{T}} {\ast}_{M} \dfrac{\partial \bm{\mathcal{X}}^{T}}{\partial \bm{\mathcal{V}}^{N}} = \bm{\mathcal{U}}^{N} {\ast}_{M} (\delta \bm{\mathcal{X}})^{T}
\end{equation}

\subsubsection{Derivation for Multi-Layer Perception}


Using Definition \ref{*_{M}-product gradient}, we can write the back-propagation formulas for multi-layer perception. 
(Latent Horizontal Slice)
Let \bm{$\vec{\mathcal{Z}}}_{i}^{j}$ = \bm{$\mathcal{W}}_{u}^{j}$ $\ast_{M}$ \bm{$\vec{\mathcal{U}}}_{i}^{j}$ + \bm{$\vec{\mathcal{B}}}_{u}^{j}$
or 
$\mathcal{Z}_{i}^{j} = \bm{\mathcal{W}}_{u}^{j}  \ast_{M}$ \bm{$\mathcal{U}}^{j}$ + \bm{$\vec{\mathcal{B}}}_{u}^{j}$, the non-activated features. Then,

\begin{equation}\label{key_delta_U^j}
	\delta \bm{\mathcal{U}}^{j} = \dfrac{\partial E}{\partial \bm{\mathcal{U}}^{j}} = \dfrac{\partial E}{\partial \bm{\mathcal{U}}^{j+1}}
	\ast_{M} \dfrac{\partial \bm{\mathcal{U}}^{j+1}}{\partial \bm{\mathcal{U}}^{j}} = \delta \bm{\mathcal{U}}^{j+1} \ast_{M} 
	\dfrac{\partial \bm{\mathcal{U}}^{j+1}}{\partial \bm{\mathcal{U}}^{j}} =
	(\bm{\mathcal{W}}_{u}^{j})^{T} \ast_{M} [\sigma^{'j}(\bm{\mathcal{Z}}_{j})
	\odot \delta \bm{\mathcal{U}}^{j+1}]
\end{equation}

\begin{equation}\label{key15}
	\delta \bm{\mathcal{W}}_{u}^{j} = \dfrac{\partial E}{\partial \bm{\mathcal{W}}_{u}^{j}} = \dfrac{\partial E}{\partial \bm{\mathcal{U}}^{j+1}}
	\ast_{M} \dfrac{\partial \bm{\mathcal{U}}^{j+1}}{\partial \bm{\mathcal{W}}_{u}^{j}} = \delta \bm{\mathcal{U}}^{j+1} \ast_{M} 
	\dfrac{\partial \bm{\mathcal{U}}^{j+1}}{\partial \bm{\mathcal{W}}_{u}^{j}} =
	  [\sigma^{'j}(\bm{\mathcal{Z}}_{j}) \odot \delta \bm{\mathcal{U}}^{j+1}] \ast_{M} (\bm{\mathcal{U}}^{j})^{T}
\end{equation}

\begin{equation}\label{key16}
	\delta \bm{\vec{\mathcal{B}}}_{u}^{j} = \dfrac{\partial E}{\partial \bm{\vec{\mathcal{B}}}_{u}^{j}} = \dfrac{\partial E}{\partial \bm{\mathcal{U}}^{j+1}}
	\ast_{M} \dfrac{\partial \bm{\mathcal{U}}^{j+1}}{\partial \bm{\vec{\mathcal{B}}}_{u}^{j}} = \delta \bm{\mathcal{U}}^{j+1} \ast_{M} 
	\dfrac{\partial \bm{\mathcal{U}}^{j+1}}{\partial \bm{\vec{\mathcal{B}}}_{u}^{j}} =
	{\rm sum}(\sigma^{'j}(\bm{\mathcal{Z}}_{j}) \odot \delta \bm{\mathcal{U}}^{j+1}, 2)
\end{equation}

for $j=1,\cdots,N$ where $\odot$ is the Hadamard element-wise product and $\sigma'_{j}$
is the derivative of the activation function, applied element-wise. The function ${\rm sum}(\cdot,2)$ sums
the tensor along the second dimension. These formulas are analogous to the matrix
formulas presented in \cite{N2017}, and further highlight the benefits of matrix-mimeticity.

(Latent Lateral Slice)
Let \bm{$\vec{\mathcal{Z}}}_{i}^{j}$ = \bm{$\mathcal{W}}_{v}^{j}$ $\ast_{M}$ \bm{$\vec{\mathcal{V}}}_{i}^{j}$ + \bm{$\vec{\mathcal{B}}}_{v}^{j}$
or 
\bm{$\mathcal{Z}}_{i}^{j}$ = \bm{$\mathcal{W}}_{v}^{j}$ $\ast_{M}$ \bm{$\mathcal{V}}^{j}$ + \bm{$\vec{\mathcal{B}}}_{v}^{j}$ , the non-activated features. Then,

\begin{equation}\label{key_delta_V^j}
	\delta \bm{\mathcal{V}}^{j} = \dfrac{\partial E}{\partial \bm{\mathcal{V}}^{j}} = \dfrac{\partial E}{\partial \bm{\mathcal{V}}^{j+1}}
	\ast_{M} \dfrac{\partial \bm{\mathcal{V}}^{j+1}}{\partial \bm{\mathcal{V}}^{j}} = \delta \bm{\mathcal{V}}^{j+1} \ast_{M} 
	\dfrac{\partial \bm{\mathcal{V}}^{j+1}}{\partial \bm{\mathcal{V}}^{j}} =
	(\bm{\mathcal{W}}_{v}^{j})^{T} \ast_{M} [\sigma^{'j}(\bm{\mathcal{Z}}_{j})
	\odot \delta \bm{\mathcal{V}}^{j+1}]
\end{equation}

\begin{equation}\label{key18}
	\delta \bm{\mathcal{W}}_{v}^{j} = \dfrac{\partial E}{\partial \bm{\mathcal{W}}_{v}^{j}} = \dfrac{\partial E}{\partial \bm{\mathcal{V}}^{j+1}}
	\ast_{M} \dfrac{\partial \bm{\mathcal{V}}^{j+1}}{\partial \bm{\mathcal{W}}_{v}^{j}} = \delta \bm{\mathcal{V}}^{j+1} \ast_{M} 
	\dfrac{\partial \bm{\mathcal{V}}^{j+1}}{\partial \bm{\mathcal{W}}_{v}^{j}} =
	[\sigma^{'j}(\bm{\mathcal{Z}}_{j}) \odot \delta \bm{\mathcal{V}}^{j+1}] \ast_{M} (\bm{\mathcal{V}}^{j})^{T}
\end{equation}

\begin{equation}\label{key19}
	\delta \bm{\vec{\mathcal{B}}}_{v}^{j} = \dfrac{\partial E}{\partial \bm{\vec{\mathcal{B}}}_{v}^{j}} = \dfrac{\partial E}{\partial \bm{\mathcal{V}}^{j+1}}
	\ast_{M} \dfrac{\partial \bm{\mathcal{V}}^{j+1}}{\partial \bm{\vec{\mathcal{B}}}_{v}^{j}} = \delta \bm{\mathcal{V}}^{j+1} \ast_{M} 
	\dfrac{\partial \bm{\mathcal{V}}^{j+1}}{\partial \bm{\vec{\mathcal{B}}}_{v}^{j}} =
	{\rm sum}(\sigma^{'j}(\bm{\mathcal{Z}}_{j}) \odot \delta \bm{\mathcal{V}}^{j+1}, 2)
\end{equation}

\section{\large\bf Experiments}

In this section, BDTFN is tested using video and color images. The RMSE value and running time are compared using BDTFN with three state-of-the art algorithms such as TNN\cite{ZA2017}, TLNM-TQR\cite{Yongming2016} and LATC-Tubal\cite{ChenXY2020}. We also evaluate the data completion performance of the BDTFN model using two large-scale and high-dimensional traffific data sets collected from California PeMS system. We measure our model by both imputation accuracy and computational cost.

All our methods have been implemented for tensor completion on the Pycharm 2018.2.5 x64 platform equipped with an Intel Core i5-8300H CPU, 8GB of RAM and python 3.7.4.

\subsection{California PeMS Data Sets}
To show the advantages of BDTFN for handling with large-scale and high-dimensional traffific data, we particularly choose two publicly available data sets collected from California transportation system (i.e., PeMS) as our benchmark data sets\footnote{
	The data sets are available at \tcb{\url{https://doi.org/10.5281/zenodo.3939792}}}

\begin{itemize}
	\item PeMS-4W data set: This data set contains freeway traffific speed collected from 11160 traffic measurement sensors over 4 weeks (the fifirst 4 weeks in the year of 2018) with a 5-minute time resolution (288 times intervals per day) in California, USA. It can be arranged in a matrix of size 11160 $\times$ 8064 or a tensor of size 11160 $\times$ 288 $\times$ 28 according to the spatial and temporal dimensions. Note that this data set contains about 90 million observations. 
	\item PeMS-8W data set: This data set contains freeway traffific speed collected from 11160 traffic measurement sensors over 8 weeks (the fifirst 8 weeks in the year of 2018) with a 5-minute time resolution (288 times intervals per day) in California, USA. It can be arranged in a matrix of size 11160 $\times$ 16128 or a tensor of size 11160 $\times$ 288 $\times$ 56 according to the spatial and temporal dimensions. Note that this data set contains about 180 million observations.
\end{itemize}

It is not diffificult to see that these two data sets are both large-scale and high-dimensional. In what follows, we create a missing patterns, i.e., random missing (RM), which are same as the work \cite{ChenSun2020}. Then according to the mechanism of RM patterns, we mask certain amount of observations as missing values (i.e., 30$\%$, 70$\%$) in both two data sets, and the remaining partial observations are input data for imputing these masked entries. To assess the imputation performance, we use the actual values of the masked missing entries as the ground truth to compute the matrices RMSE: 

\begin{equation}
	RMSE = \sqrt{\dfrac{1}{n} \sum\limits_{i=1}^{n} (y_{i} - \hat{y}_{i})^{2}},
\end{equation}

where $y_{i}, \hat{y}_{i}, i=1,...,n$ are actual values and estimated/imputed value.

Table \ref{Performance comparison (in RMSE) for RM data imputation tasks on California PeMS-4W and PeMS-8W traffic speed data. } shows the results of BDTFN and its competing models for missing traffic data imputation tasks. For RM scenarios on both PeMS-4W and PeMS-8W data sets, the proposed BDTFN model achieves high accuracy. By comparing with LATC-Tubal, we can see that BDTFN performs consistently better in all missing scenarios. This clearly shows that tensor decomposition based on deep learning is superior to weighted tensor kernel norm minimization.

\begin{table}[H]
	\caption{Performance comparison (in RMSE) for RM data imputation tasks on California PeMS-4W and PeMS-8W traffic speed data. }
	\label{Performance comparison (in RMSE) for RM data imputation tasks on California PeMS-4W and PeMS-8W traffic speed data. }
	\begin{center} 
		\begin{tabular}{l|ll|ll}
			\toprule
			& \multicolumn{2}{c|}{PeMS-4W} & \multicolumn{2}{l}{PeMS-8W} \\
			Completion Approach & 30 $\%$ & 70 $\%$ & 30 $\%$ & 70 $\%$  \\
			\midrule
			TNN                    & 63.52          & 63.51          & 63.58          & 63.57 \\
			TLNM-TQR                      & 63.51          & 63.50          & 63.57          & 63.57 \\
			LATC-Tubal                    & 1.71           & 2.46           & 1.74           & 2.51  \\
			BDTFN                          & \textbf{1.07}  & \textbf{1.14}  & \textbf{1.05}  & \textbf{0.84} \\
			\bottomrule
		\end{tabular}
	\end{center}
\end{table}

In Table \ref{Performance comparison (in RMSE) for RM data imputation tasks on California PeMS-4W and PeMS-8W traffic speed data. }, we see that BDTFN and LATC-Tubal can produce comparable imputation accuracy. However, observing the running time of imputation models in Table \ref{Running time  (in minute)  of tensor completion result for California PeMS-4W and PeMS-8W traffic speed data.}, it is not difficult to conclude that both TNN and TLNM-TQR are not suitable for these large-scale imputation tasks due to their high computational cost. Table \ref{Running time  (in minute)  of tensor completion result for California PeMS-4W and PeMS-8W traffic speed data.} clearly shows that BDTFN with tensor decomposition based on deep learning is the most computationally efficient one by comparing to the competing models. Therefore, from both Table \ref{Performance comparison (in RMSE) for RM data imputation tasks on California PeMS-4W and PeMS-8W traffic speed data. } and \ref{Running time  (in minute)  of tensor completion result for California PeMS-4W and PeMS-8W traffic speed data.}, the results suggest that the proposed BDTFN model is an efficient solution to large-scale traffic data imputation while still maintaining the promising imputation performance close to the state-of-the-art models.

\begin{table}[H]
	\caption{Running Time (in minute)  of tensor completion result for California PeMS-4W and PeMS-8W traffic speed data.}
	\label{Running time  (in minute)  of tensor completion result for California PeMS-4W and PeMS-8W traffic speed data.}
	\begin{center} 
		\begin{tabular}{l|ll|ll}
			\toprule
			 & \multicolumn{2}{c|}{PeMS-4W} & \multicolumn{2}{l}{PeMS-8W} \\
			 & 30 $\%$ & 70 $\%$ & 30 $\%$ & 70 $\%$  \\
			\midrule
			TNN                    & 82.58(20)        & 81.35(20)        & 169.39(20)         & 165.78(20)    \\
			TLNM-TQR                      & 84.23(20)        & 82.15(20)        & 165.33(20)         & 167.70(20)      \\
			LATC-Tubal                    & 16.96(48)       & 29.10(69)       & 42.75(48)         & 74.77(70) \\
			BDTFN         & \textbf{12.35}(200)   & \textbf{12.65}(200)   & \textbf{28.65}(200)   & \textbf{26.81}(200)  \\
			\bottomrule
		\end{tabular}
	\end{center}
\end{table}

\begin{figure}[H]
	\centering
	\subfigure[The 1st time series.]{
		\label{Fig.sub.1}
		\includegraphics[width=1\textwidth]{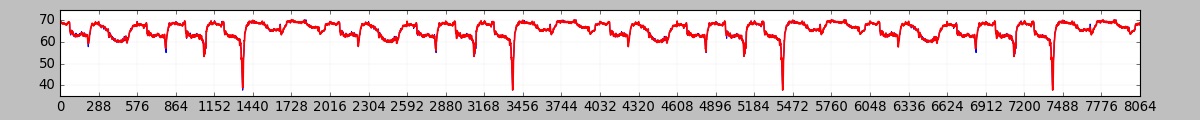}}\\
	\subfigure[The 100th time series.]{
		\label{Fig.sub.2}
		\includegraphics[width=1\textwidth]{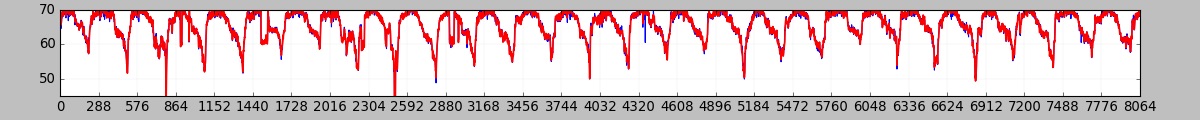}}\\
	\subfigure[The 200th time series.]{
		\label{Fig.sub.3}
		\includegraphics[width=1\textwidth]{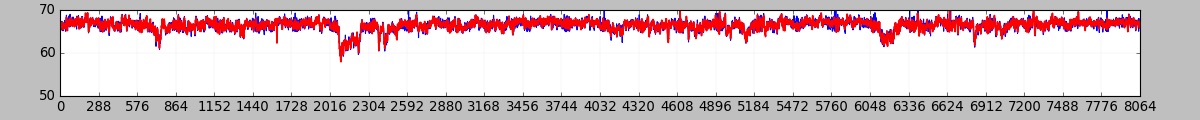}}
	\caption{Estimated time series (red curve) vs. ground truth (blue curve) for the PeMS-4W data set at the case of 30$\%$ RM scenario.}
	\label{time series}
\end{figure}

Figure \ref{time series} shows that BDTFN can produce masked time series points accurately by learning from partial observations.

\subsection{Color Image Recovery}

In this part, we use a $n_{1}\times n_{2}$ size color image as a $n_{1}\times n_{2}\times 3$ tensor to test the performance of our algorithm and compare it with other state-of-the-art algorithms with the same parameters as were previously set. All of the color images are from the Berkeley Segmentation Dataset \cite{MFTM2001} and their sizes are $481\times 321$ or $321\times 481$. 

\begin{figure}[H]
	\centering
	\noindent\makebox[\textwidth][c] {
		\includegraphics[scale=0.36]{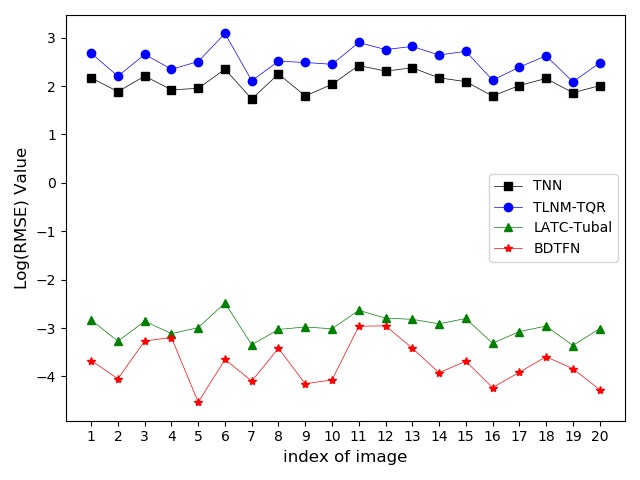} }
	\caption{Log(RMSE) Value of the recovery results for the four algorithms.}
	\label{Color Image Recovery RMSE Value}
\end{figure}

\begin{figure}[H]
	\centering
	\noindent\makebox[\textwidth][c] {
		\includegraphics[scale=0.36]{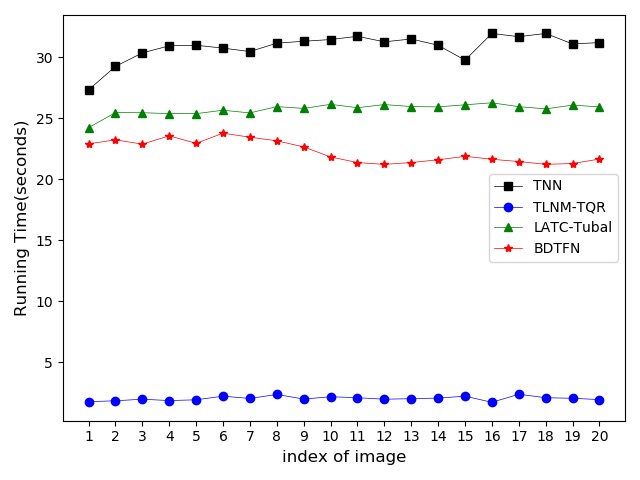} }
	\caption{Running Time of the recovery results for the four algorithms.}
	\label{Color Image Recovery Running Time}
\end{figure}

\begin{figure}[H]
	\centering
	\noindent\makebox[\textwidth][c] {
		\includegraphics[scale=0.7]{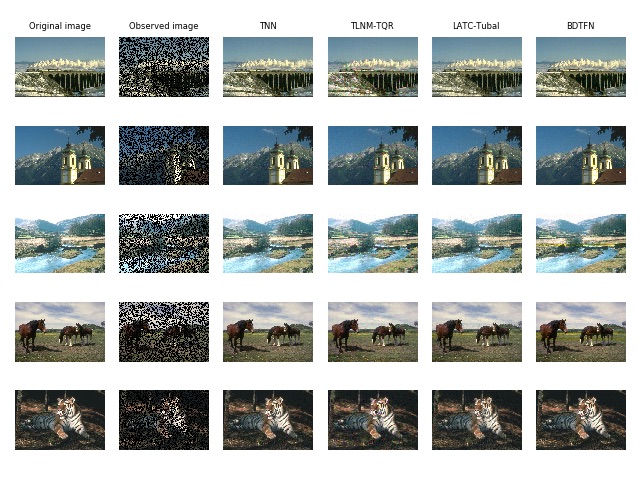} }
	\caption{Recovery performance comparison on 5 example images.}
	\label{Recovery image show}
\end{figure}

To further test the performance of our algorithm, the Python function, imnoise, is used to add a blurring noise to the 20 color images as well as add Gaussian noise with a mean of zero and a standard deviation $\sigma=5e-3$ to make the images more difficult to recover. In addition, we maintain 50\% of the given pixels as known entries in $\Omega$. The results of the tests and comparison are shown in Figure \ref{Color Image Recovery RMSE Value} and \ref{Color Image Recovery Running Time}, and some examples of the recovered images are shown in Figure \ref{Recovery image show}.

From these results and the experimental data in Table \ref{Comparison results of running time and the RMSE on the 5 images}, it can be seen that BDTFN is still the fastest and that its image recovery capabilities are satisfactory.

\begin{table}[H]
	\caption{Comparison of the results of running time and the RMSE on the 5 color images in Figure \ref{Recovery image show}}
	\label{Comparison results of running time and the RMSE on the 5 images}
	\resizebox{\textwidth}{!}
	{
		\begin{tabular}{ccccccccc}
			\toprule
			\multirow{2}{*}{Index} &  \multicolumn{4}{c}{Running Time (seconds)}&\multicolumn{4}{c}{RMSE Value}\\  
			\cmidrule(lr){2-5} \cmidrule(lr){6-9}
			&TNN &TLNM-TQR  &LATC-Tubal  &BDTFN  &TNN &TLNM-TQR  &LATC-Tubal  &BDTFN \\
			\midrule
			1 & 25.7344  & \textbf{1.7500}  & 23.2000  & 21.5284  & 8.7349  & 14.5508  & 0.0593  & \textbf{0.0266} \\
			2 & 24.0469  & \textbf{1.7813}  & 22.5400  & 20.6648  & 6.5897  & 9.0614  & 0.0379  & \textbf{0.0238} \\
			3 & 27.2500  & \textbf{1.9375}  & 23.3600  & 20.0724  & 9.0877  & 14.2380  & 0.0577  & \textbf{0.0387} \\
			4 & 27.5625  & \textbf{1.9687}  & 22.9200  & 20.9201  & 6.8410  & 10.4619  & 0.0573  & \textbf{0.0368} \\
			5 & 28.7031  & \textbf{1.8906}  & 23.2400  & 20.8473  & 7.0441  & 12.2631  & 0.0504  & \textbf{0.0106} \\
			\bottomrule
	\end{tabular}}
	
\end{table}

\subsection{Video Recovery}
Here we use a black and white gray-scale basketball video, which depicts 1.6 seconds of the game and contains 40 digital images in AVI format. Each frame of this video holds $144\times256$ pixels of black and white images so that it can be considered as a three-dimensional tensor $\mathcal{X}\in \mathbb{R}^{144\times 256\times 40}$. If we retain $50\%$ of the given pixels as known entries in $\Omega$, the problem of video recovery problem will be considered a tensor completion problem. Thus, we can solve this problem using our BDTFN algorithm.

To verify that BDTFN is convergent, we chose a video with $50\%$ miss rate to test our algorithm. When all the parameters are fixed, the required inputs in BDTFN are the iteration number and tubal-rank of the result. Therefore, we respectively fixed the iteration number and tubal-rank in Figure 5. It is easy to see that our method has a good convergence accuracy and only needs a few iterations to converge, therefore showing our algorithm to be fast and precise. This also gives us confidence to transcend other algorithms that have been popular in recent years.

\begin{figure}[H]
	\centering
	\includegraphics[scale=0.7]{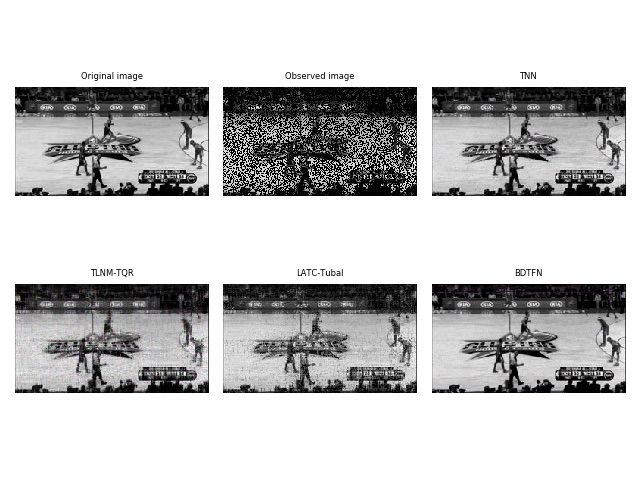} 
	\caption{The 30th frame of the completion result in a basketball video.}
	\label{The 30th frame of completion result for a basketball video}
\end{figure}	

In order to achieve the optimal effect of the algorithms, we set the parameter $\lambda$ of TNN to $\lambda =\frac{1}{\sqrt{3\max \left ( n_{1},n_{2} \right ) } } $ and adjust its iteration to 20. For TLNM-TQR, we set $r=11$ and its parameter $\mu$ and $\rho$ are adjusted to be $\mu=10^{-2}$ and $\rho=1.5$ respectively. From Figure \ref{The 30th frame of completion result for a basketball video}, our algorithm can be seen to be stable while the miss rate is $50\%$. The detail comparison data are listed in Table \ref{Running time and RMSE of tensor completion result for the basketball video}.

\begin{table}[H]
	\caption{Running time and RMSE of tensor completion result for the basketball video}
	\label{Running time and RMSE of tensor completion result for the basketball video}
	\begin{center} 
		\begin{tabular}{lcc}
			\toprule
			Completion Approach & RMSE Value &Running Time (seconds)    \\
			\midrule
			TNN      & 8.1824           & 295.6094 \\
			TLNM-TQR        & 12.7794          & \textbf{34.2188} \\
			LATC-Tubal      & 0.0911           & 233.9589 \\
			BDTFN            & \textbf{0.0219}  & 106.9024 \\
			\bottomrule
		\end{tabular}
	\end{center}
\end{table}

From these results, we can see that the speed of our algorithm is slower than TLNM-TQR, while its accuracy is faster then others, which can also be seen in Figure \ref{The 30th frame of completion result for a basketball video}. Further to this, we explore the performance of our algorithm using a variety of specified missing ratio from $30\%$ to $90\%$. All the algorithms are tested 50 times to reduce contingency.

\begin{figure}[H]
	\centering
	\subfigure{
		\includegraphics[width=0.5\textwidth]{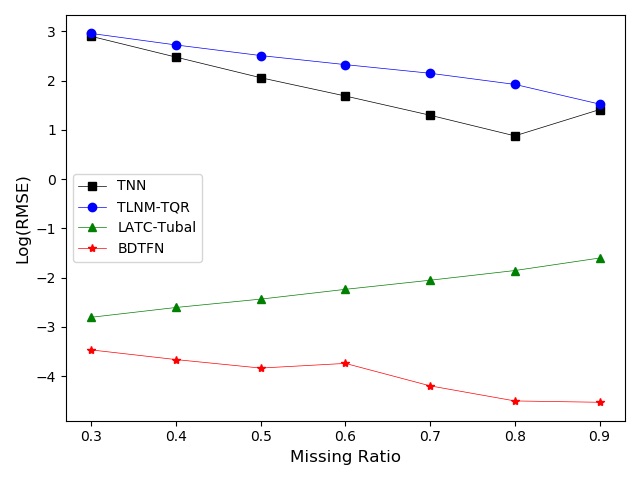}}
	\caption{Runnning Time of the four algorithms with different missing ratios.}
	\label{RMSE Value comparative result of BDTFN}
\end{figure}

\begin{figure}[H]
	\centering
	\subfigure{
		\includegraphics[width=0.5\textwidth]{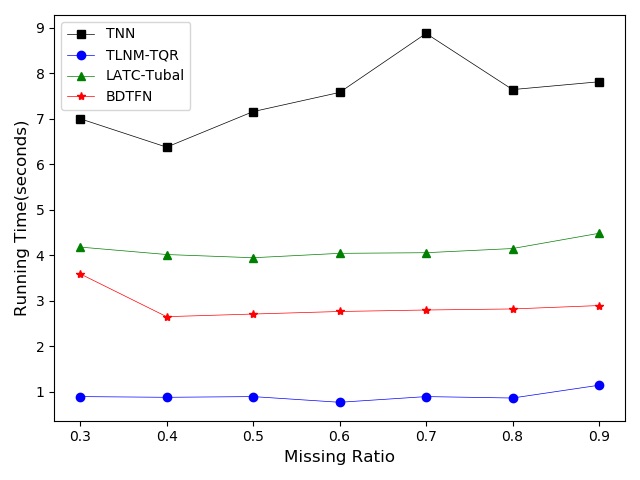}}
	\caption{Runnning Time of the four algorithms with different missing ratios.}
	\label{Running Time comparative result of BDTFN}
\end{figure}

As can be seen from Figure \ref{RMSE Value comparative result of BDTFN} and \ref{Running Time comparative result of BDTFN}, LATC-Tubal is as precise as BDTFN and more precise than others for each specified missing ratio. 

\begin{figure}[H]
	\centering
	\subfigure{ 
		\includegraphics[width=0.5\textwidth]{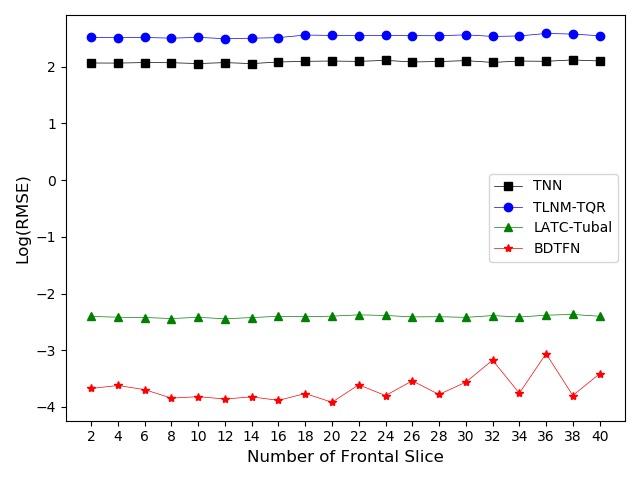} }
	\caption{RMSE Value of the four algorithms with different numbers of frontal slice.} 
	\label{RMSE Value frontal comparative result of BDTFN}
\end{figure}

\begin{figure}[H]
	\centering
	\subfigure{ 
		\includegraphics[width=0.5\textwidth]{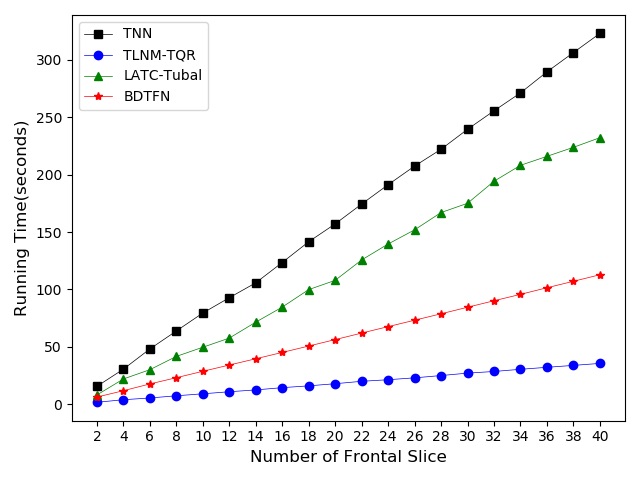} }
	\caption{RMSE Value of the four algorithms with different numbers of frontal slice.} 
	\label{Running Time frontal comparative result of BDTFN}
\end{figure}

For the next experiment, we are going to set a variety of frontal slice $\left( \text{video frame} \right)$ numbers from 2 to 40, in order to provide a recovery performance comparison of the four algorithms with the same experiment setup as the above-mentioned. As shown in Figure \ref{RMSE Value frontal comparative result of BDTFN} and \ref{Running Time frontal comparative result of BDTFN}, when the number of frontal slices steadily increases, the RMSE result of BDTFN is close to zero and better than other algorithms. In addition, the running time line graph shows that the rising rate of BDTFN is second only to TLNM-TQR, which shows that our algorithm will have great advantages in dealing with more complex tensor structures and data.

\section{Conclusion}
In this work, we propose efficient and scalable tensor completion and tensor sensing algorithms. In order to be suitable for large-scale data and obtain high convergence, we define a new deep neural network model. The experimental results show that the model has high efficiency and high speed in large-scale space-time interpolation. BDTFN is very efficient in video and color image restoration.

\end{document}